# A Multi-scale Perimeter Control and Route Guidance System for Large-scale Road Networks

Xianyue Peng, Hao Wang, Shenyang Chen, H. Michael Zhang

*Abstract*— **Perimeter control and route guidance are effective ways to reduce traffic congestion and improve traffic efficiency by controlling the spatial and temporal traffic distribution on the network. This paper presents a multi-scale joint perimeter control and route guidance (MSJC) framework for controlling traffic in large-scale networks. The network is first partitioned into several subnetworks (regions) with traffic in each region governed by its macroscopic fundamental diagram (MFD), which forms the macroscale network (upper level). Each subnetwork, comprised of actual road links and signalized intersections, forms the microscale network (lower level). At the upper level, a joint perimeter control and route guidance model solves the region-based inflow rate and hyper-path flows to control the accumulation of each region and thus maximize the throughput of each region. At the lower level, a perimeter control strategy integrated with a backpressure policy determines the optimal signal phases of the intersections at the regional boundary. At the same time, a route choice model for vehicles is constructed to meet hyper-path flows and ensure the intra-region homogeneity of traffic density. The case study results demonstrate that the proposed MSJC outperforms other benchmarks in regulating regional accumulation, thereby improving throughput.**

## I. INTRODUCTION

Perimeter control and route guidance are effective measures to control traffic flow and avoid hyper-congestion. Considering that they both affect the traffic distribution and congestion propagation in the road network, it is worthwhile to examine the relationship between perimeter control and route guidance and develop methods to jointly optimize them to improve the performance of large-scale traffic networks. Given the high complexity involved in controlling these networks, a common approach is to divide them into several homogeneous subnetworks (regions). Macroscopic fundamental diagram (MFD) has become a popular model used in developing traffic controls at the regional level, as it provides a unimodal, low-scatter and demand-insensitive relationship between vehicle accumulation and outflow [1, 2]. Consequently, there have been numerous studies focusing on MFD-based perimeter control and MFD-based route guidance.

To help review the relevant work in this area, we provide a general framework for multi-scale perimeter control and route guidance system based on MFD. This control system consists of two spatial/temporal scales:

· *Macroscopic scale (upper level)*: employs perimeter control and hyper-path1 flow assignment for the entire network.

· *Microscopic scale (lower level)*: utilizes signal control at regional boundaries and carries out vehicle routing within a region.

At the upper level, perimeter control regulates the overall inflow of vehicles into a region, while the lower level focuses on modeling and controlling specific entries, such as intersections at the regional boundary. Route guidance at the upper level determines the flows between regions, known as hyper-path flows, while vehicle routes consisting of actual road links are determined at the lower level.

In the context of the above framework, existing studies often fall in one of the following three categories: one- or two-level perimeter control; one- or two-level route guidance; joint perimeter control and route guidance, mainly at the upper level. For example, many MFD-based perimeter control methods optimize the inflow rate at the macroscopic scale to control accumulation in each region of the network [3-5]. As for the one-region system, Keyvan-Ekbatani, et al. [3] designed a real-time gating control method using feedback control, Haddad and Shraiber [4] proposed a robust controller to cope with uncertainties in traffic. Additionally, de Souza, et al. [6] provided a comprehensive comparison between feedback perimeter control and traffic-responsive urban control, highlighting how their combination improves network performance. Two-region and multi-region perimeter controls

This work was supported by the National Key R&D Program of China [2022ZD0115600] and National Natural Science Foundation of China [No. 52472341] (*Corresponding author: Hao Wang*).

Xianyue Peng is with the School of Transportation, Southeast University, Nanjing 211189, China and also with the Department of Civil and Environmental Engineering, University of California, Davis, CA 95616, USA (e-mail: xypeng@seu.edu.cn)

Hao Wang is with Jiangsu Key Laboratory of Urban ITS, Southeast University, Nanjing 211189, China; Jiangsu Province Collaborative Innovation Center of Modern Urban Traffic Technologies, Nanjing 211189, China; School of Transportation, Southeast University, Nanjing 211189, China (e-mail: haowang@seu.edu.cn).

Shenyang Chen is with the Department of Civil and Environmental Engineering, University of California, Davis, CA 95616, USA (e-mail: aschen@ucdavis.edu).

H. Michael Zhang is with the Department of Civil and Environmental Engineering, University of California, Davis, CA 95616, USA (e-mail: hmzhang@ucdavis.edu).

[1] We use the term hyper-path rather than route for a "route" at the upper level to distinguish it with a route at the lower level, since the former comprises a sequence of regions rather than actual road links, while the latter comprises actual road links.

have also been investigated. Geroliminis, et al. [5] developed an optimal perimeter control method for a two-region system to maximize the trip completion. Jiang and Keyvan-Ekbatani [7] integrated a real-time partitioning into the perimeter control system and introduced a Deep Learning estimator to reduce computational burden. In addition, Haddad and Mirkin [8], Ampountolas, et al. [9], Zhong, et al. [10], Wang, et al. [11], Bichiou, et al. [12], Ding, et al. [13] proposed adaptive perimeter control, bi-modal perimeter control, robust perimeter control, perimeter control for the road networks with incidents, sliding mode perimeter control, and perimeter control based on three-dimensional MFD respectively.

An important issue that has not been adequately addressed in these studies is how to map the perimeter control variables to signal timings at the boundary intersections to achieve the desired inflow rate. Aboudolas and Geroliminis [14] provided a calculation equation to determine the green split at boundary intersections based on the saturated flow rate. Mousavizadeh and Keyvan-Ekbatani [15] proposed a method to estimate regional transfer flow ratios based on floating car data. Keyvan-Ekbatani, et al. [16] addressed the delay and queue management for the gated links at the boundary intersection using a quadratic knapsack model. Yang, et al. [17] introduced a multi-scale perimeter control model that minimizes travel costs at both local and network levels by optimizing green split at boundary intersections. Tsitsokas, et al. [18] adopted proportional–integral regulator for MFD-based gating and formulated an optimization problem to calculate the green time. An alternative approach involved a hierarchical approach with a high-level controller based on a regional model and a low-level controller based on a sub-regional model [19, 20].

MFD-based route guidance assigns traffic demand to each hyper-path to achieve a specific objective, such as system optimal or user equilibrium. Leclercq and Geroliminis [21], Gayah and Daganzo [22] analyzed the effect of route choice on MFD in simple networks. Hajiahmadi, et al. [23] proposed a dynamic route guidance method to obtain the hyper-path flows for a multi-region network. Menelaou, et al. [24], Menelaou, et al. [25] integrated multi-regional route guidance with demand management. To obtain more detailed routes for vehicles, Yildirimoglu and Geroliminis [26], Yildirimoglu, et al. [27] integrated dynamic traffic assignment including user equilibrium and system optimum to the MFD-based model, and then proposed a hierarchical structure comprising region-level and sub-region-level MFD-based models to improve network performance and provide sub-regional vehicle paths, respectively[28]. Moreover, Jiang, et al. [29] proposed a regional route guidance framework based on model predictive control and deep reinforcement learning to generate customized route plans while accounting for drivers' compliance.

Recent research has focused on integrating perimeter control with route guidance at the upper level to optimize both inflow rate and route proportion, thus improving overall network performance. Various approaches have been proposed to achieve this integration, including MFD-based models proposed by Sirmatel and Geroliminis [30] and data-driven model-free control strategies proposed by Hou and Lei [31]. Ding, et al. [1] developed separate algorithms for traffic guidance (based on variable message signs) and perimeter control, which are sequentially executed until the trip completion flow is maximized. Similar to the perimeter control structure proposed by Ramezani, et al. [19], Fu, et al. [20], the hierarchical approach was also employed in the joint control strategy [32, 33]. Fu, et al. [32] optimized the inflow rate for a two-region model at the upper level and then solved for the inflow rate and route proportion at the lower subregion level. A colored petri nets-based model was proposed to represent traffic dynamics, with buffers used to model the regional boundaries and their storage capacities. Chen, et al. [33] developed an iterative adaptive dynamic programming (IADP) approach to address the model-plant mismatch, where the region-based model represents the model, and the detailed subregion model serves as the plant. Additionally, several studies have proposed perimeter control methods that incorporate route choice models to capture dynamic traffic behavior, such as those by Ramezani, et al. [19], Guo and Ban [34] and Lei, et al. [35]. Most of these studies solve the combined perimeter control and route guidance problem at the upper level, and the solution cannot be directly applied to provide routes to traffic at the lower level (within a region). Specifically, many strategies assign the flow with region-based ODs to a hyper-path rather than an actual route that vehicles can travel on. Additionally, their perimeter control always focuses on the general traffic movements between adjacent regions without considering their specific movements at the boundary intersections, nor do they consider the capacities of the downstream links of the boundary intersections.

This paper aims to develop a multi-scale perimeter control and route guidance system for multi-region road networks to fill the gaps revealed in the above review. The multi-scale approach integrates the broad, long-term strategies at the macroscopic level with the detailed, short-term execution of routing and signal control at the microscopic level. At the macroscopic level, it manages overall network flow and prevents widespread congestion, while at the microscopic level, it adjusts to immediate local traffic conditions with precise control measures. This synergy ensures that comprehensive traffic strategies are effectively turned into specific, actionable steps on the ground.

Specifically, the paper introduces a joint perimeter control and route guidance model that optimizes perimeter control and route guidance control variables simultaneously; then a perimeter control algorithm combined with a multi-phase and backpressure control policy is proposed to control the boundary intersections with the help of a more accurate method for estimating the boundary throughput; and finally a route choice model is developed to obtain actual routes of vehicles that are consistent with hyper-path flow recommendations and ensures traffic homogeneity inside each region.

The remainder of this paper is organized as follows. Section II introduces the multi-scale framework for joint perimeter control and regional route guidance. In Section III, a joint control model for multi-region networks is presented. And Section IV proposes the perimeter control approach and route choice model at the microscopic scale. The case study is presented in Section V. Finally, Section VI discusses and concludes the paper.

Table I provides a summary of the key notations used throughout this paper.

TABLE I
LIST OF KEY NOTATIONS

| Category | Notation | Definition |
|---|---|---|
| **Control step** | $t, k$ | macro and micro control step respectively |
| | $K_t$ | set of micro control steps within macro control step $t$ |
| | $T_{mac}, T_{mic}$ | duration of the macro and the micro control step, respectively, s |
| **Multi-region network model** | $R$ | set of regions |
| | $R_i$ | set of the adjacent regions to region $i$ |
| | $N_i(t)$ | number of vehicles in region $i$ at macro control step $t$, veh |
| | $N_{ij}(t)$ | number of vehicles in region $i$ with destination region $j$ at macro control step $t$, veh |
| | $Q_{ij}(t)$ | newly generated traffic flow of region $i$ with destination region $j$ at macro control step $t$, veh |
| | $N_{ihj}(t)$ | number of vehicles in region $i$ with destination region $j$ which cross the boundary between region $i$ and region $h$ at macro control step $t$, , veh per micro-step |
| | $N_i^{II}(t)$ | number of *Type II vehicles* in region $i$ at macro control step $t$, veh |
| | $N_{ij}^{I}(t)$ | number of *Type I vehicles* in region $i$ with destination in region $j$ at macro control step $t$, veh |
| | $N_i^{crit}$ | critical value of region $i$ corresponding to the maximum trip completion flow, veh |
| | $G_i(N_i(t))$ | trip completion flow of vehicles in region $i$ per unit time, veh/s |
| **Joint control variables** | $b_{ih}(t)$ | perimeter control variable, representing the crossing fraction of *Type I vehicles* at the boundary between region $i$ and $h$ at macro control step $t$, % |
| | $c_{ihj}(t)$ | route guidance variable, representing the proportion of vehicles in region $i$ to destination region $j$ with region $h$ as next region to travel at macro control step $t$, % |
| | $c_{ihj}^{min}(t), c_{ihj}^{max}(t)$ | Minimum and maximum of route guidance variable, % |
| **Perimeter control from region $i$ to $h$** | $M_{ih}(t)$ | expected traffic flow rate at macro control step $t$, veh/s |
| | $m_{ih}(k)$ | expected traffic flow rate at micro control step $k$, veh/s |
| | $\widetilde{m}_{ih}(z)$ | observed actual flow rates at micro control step $z$, veh/s |
| | $M_{ih}^{min}(t), M_{ih}^{max}(t)$ | minimum and maximum expected traffic flow rate at macro control step $t$, veh/s |
| | $\widehat{m}_{ih}^{s}(k), \widehat{m}_{ih}^{s}(0)$ | estimated flow rate of the multi-phase plan $s$ at micro control step $k$ and at the first micro control step, respectively, veh/s |
| | $\widehat{m}_{ih}^{ng}(k), \widehat{m}_{ih}^{ng}(0)$ | estimated flow rate through the non-gating intersections at micro control step $k$ and at the first micro control step, respectively, veh/s |
| **Multi-phase selection** | $p$ | phase for an intersection |
| | $s$ | multi-phase for a regional boundary |
| | $P_s$ | set of phases for each intersection under the multi-phase plan $s$ |
| | $S_{ih}$ | set of multi-phase plan $s$ for the boundary $(i, h)$ |
| | $s_{ih}^{act}(k)$ | activated multi-phase for the boundary $(i, h)$ at micro control step $k$ |
| | $w_p(k)$ | weight of phase $p$ at micro control step $k$ |
| | $W_s(k)$ | weight of multi-phase plan $s$ at micro control step $k$ |
| **Route guidance** | $V_{ij}$ | set of vehicles in region $i$ destinated for region $j$ |
| | $U_v$ | set of routes of vehicle $v$ |
| | $G_r$ | upcoming region on route $r$ |
| | $X_i$ | set of links in region $i$ |
| | $\phi_{ijr}^{v}(k)$ | probability of vehicles $v$ selecting route $r$ at micro control step $k$, where $v \in V_{ij}$, % |
| | $d_x(k)$ | estimated link density of link $x$ at the end of $k$th time step, veh/m/lane |
| | $\bar{d}_i(k)$ | average link density of region $i$, veh/m/lane |

## II. FRAMEWORK

### A. Multi-scale Control System

The transportation network in the real world can be conceptualized as a complex system of multi-region networks. Each region within this system is composed of a connected road network with numerous intersections and road segments. Such a conceptualization not only prepares networks for distributed control but also structures control into two distinct levels: macroscopic (inter-region) and microscopic (intra-region) scales.

Macroscopic level control refers to strategic management of large-scale transportation networks, while microscopic level control implements these comprehensive strategies into precise, executable actions. Our multiscale control framework operates not only at two spatial scales, but also at two time scales: a macro and micro control time step, as shown in Fig.1. A macro control step $t$ can be divided into a sequence of micro control steps. These micro control steps represented by $k$ and collectively defined as the set $K_t$.

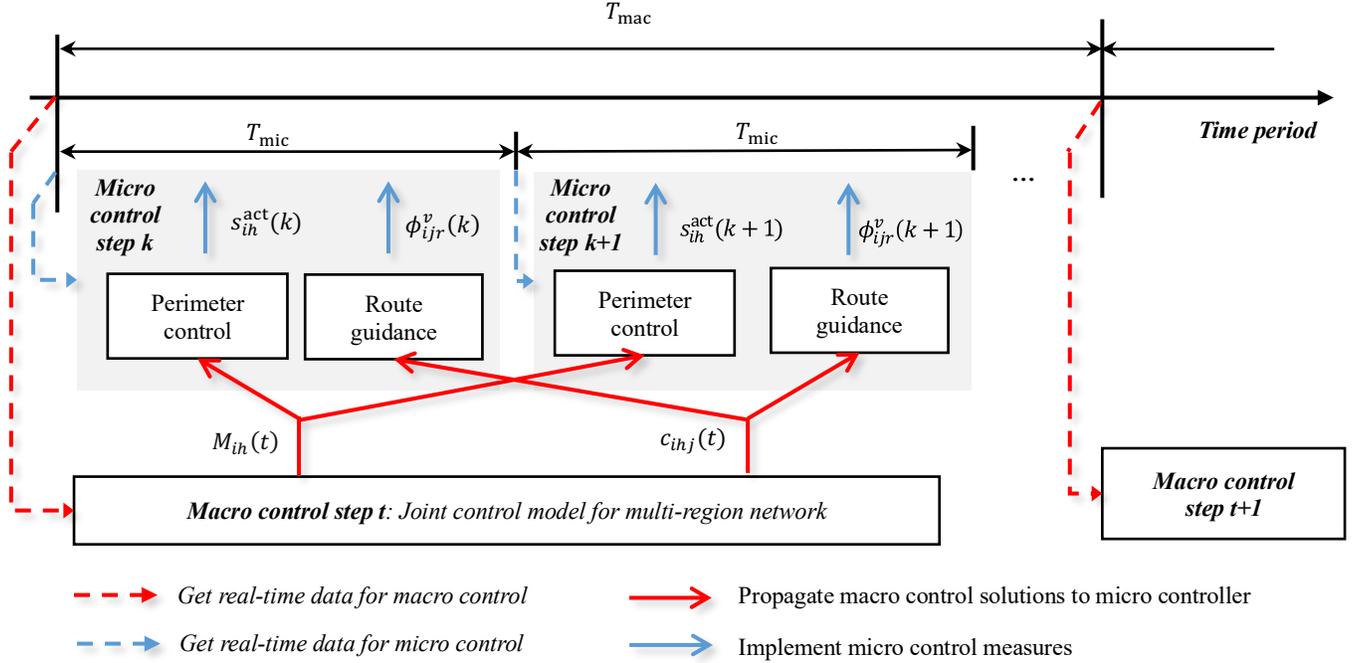

**Fig. 1.** Two temporal control scales.

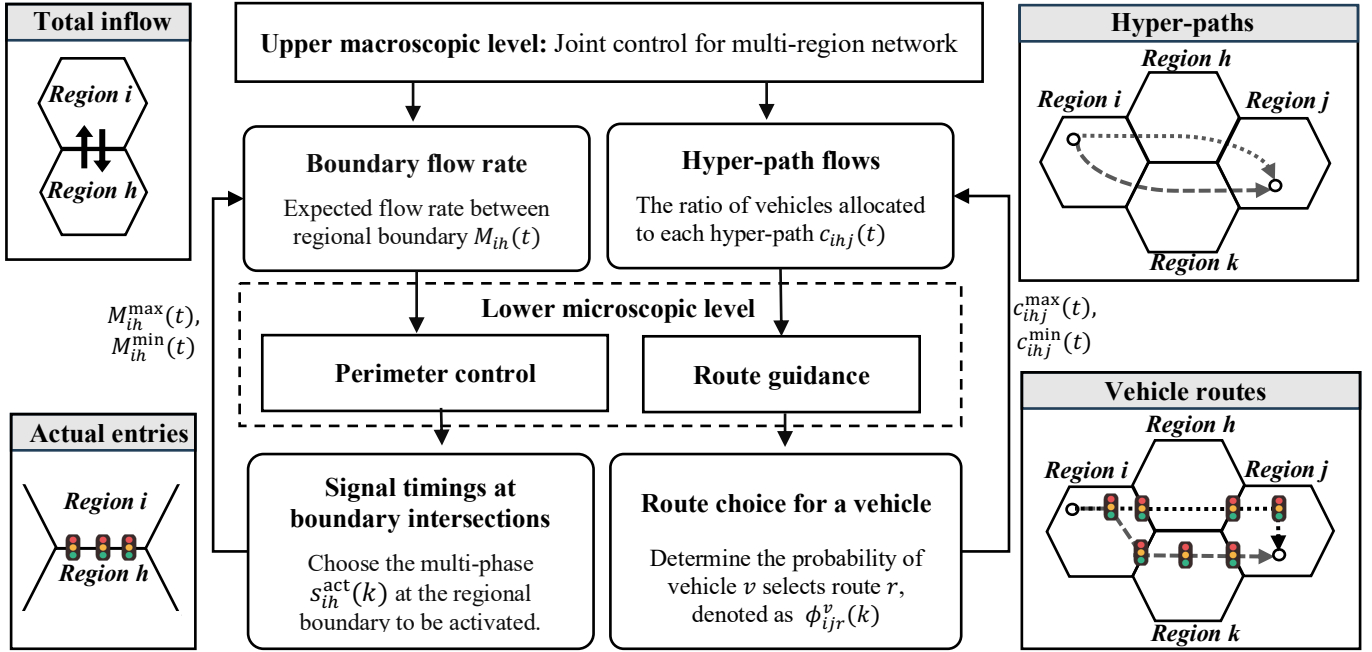

**Fig. 2.** The framework of the multi-scale control system.

The control method is based on the principles of Model Predictive Control (MPC). Specifically, at the beginning of each macro control step $t$, the joint control model, based on MFD, uses the network traffic information to predict future traffic patterns over the control horizon. This prediction allows for the optimization of control variables for the entire macro step. During the micro control step $k$, both the perimeter control model and the route choice model performs their respective control tasks.

The multi-scale perimeter control and route guidance system, illustrated in Fig. 2, consists of an upper (macroscopic) level and a lower (microscopic) level. At the upper level, a joint control model for the multi-region network manages the regional accumulation; by solving the boundary flow rate $M_{ih}(t)$ and the hyper-path flow proportions $c_{ihj}(t)$. These boundary flow rate and hyper-path flows will serve as inputs to the lower-level control. The lower-level model controls the signal timings at boundary intersections to realize the boundary flow rate from the upper level, specifically by activating the appropriate multi-phase $s_{ih}^{act}(k)$. Concurrently, it manages the route traffic inside each region to realize the hyper-path flows prescribed to it from the upper-level control, by determining the probability of vehicles selecting a route at micro control step, denoted as $\phi_{ijr}^v(k)$. Additionally, the range for the upper-level solution $M_{ih}(t)$ and $c_{ihj}(t)$ is informed by the lower-level feedback, ensuring a responsive and integrated control system.

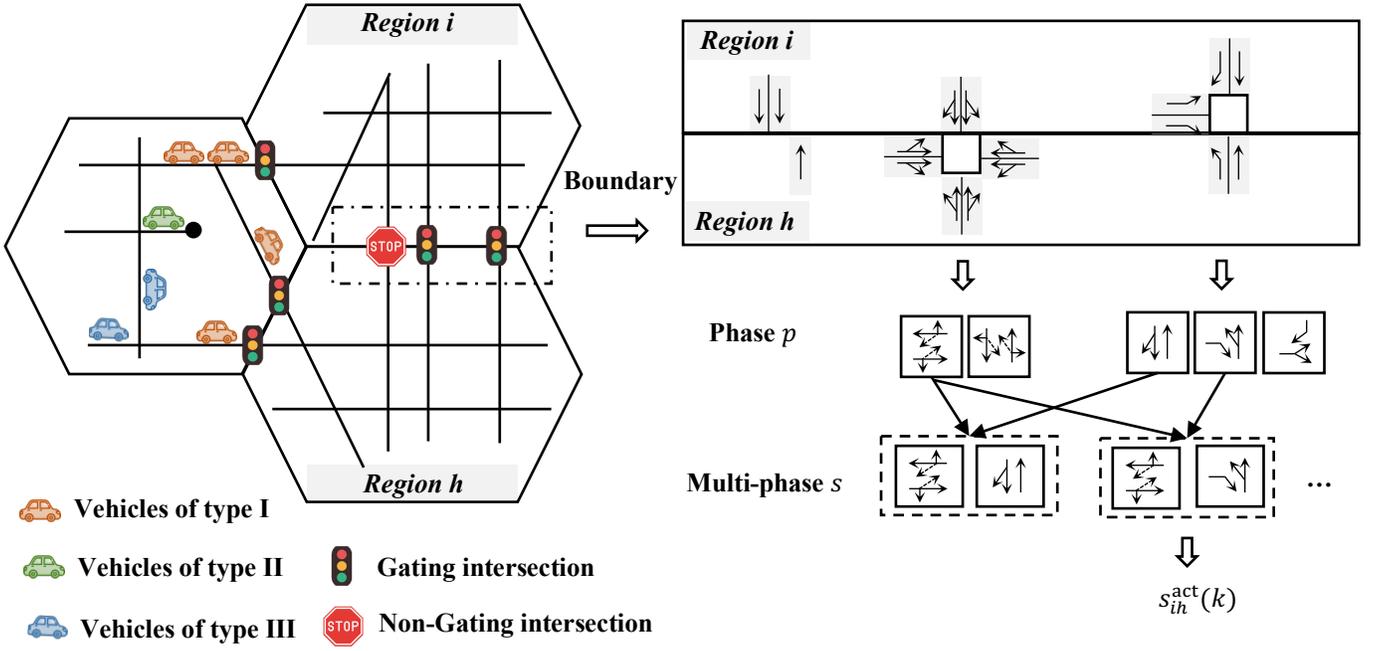

**Fig. 3.** Key components of perimeter control.

### B. Key Components

#### 1) **Perimeter Control and Multi-phase Signal Management at Regional Boundaries**

Depending on the vehicle's experience within the region, we define the vehicle states in three types:
- *Type I vehicles:* completed their journey within the region and reached regional boundary;
- *Type II vehicles:* completed their journey within the region and reached their destination;
- *Type III vehicles:* continue their movement within the region without completing their journey or reaching the region's boundary.

The number of *Type I and Type II vehicles*, which together constitutes the completion flow, can be predicted based on the accumulated number of vehicles and the MFD function. Perimeter control specifically targets *Type I vehicles*, determining whether to allow or prevent them from crossing the boundary. The joint control model at the macroscopic level is utilized to optimize the flow of these vehicles between regional boundaries.

The regional boundary comprises a series of intersections, including both gating and non-gating types. Gating intersections are controlled by the perimeter control, while non-gating intersections, such as those under priority control, are not subject to perimeter control. The flow of vehicles between each boundary solved at the macroscopic level is achieved by controlling the gating intersections at the regional boundaries at the microscopic level, specifically by determining which phases of the gating intersections are activated.

At each gating intersection, a particular phase $p$ allows certain movements and restricts others to pass through the intersection. To represent the signal plans of gating intersections at each boundary, we introduce multiple phases $s$, each corresponding to a combination of different phase $p$ settings, where $s \in S_{ih}, p \in P_s$, as illustrated in Fig.3. In practice, regional boundary control implements a multi-phase approach at each boundary during each micro control step.

#### 2) **Hyper-path Flow Redistribution and Vehicles' Rerouting**

On the other hand, rerouting is employed to redistribute network traffic at both regional and vehicle levels. For each vehicle, potential routes are initially generated using a method akin to column generation. Specifically, for each micro control step, we consider only the current route and the instantaneously optimal route for the vehicle. Following this, we create candidate hyper-paths for each vehicle based on these two routes.

At the macro level, the range of hyper-path flow is determined by aggregating the candidate hyper-paths of vehicles, as defined in Section III. The joint control model solves the hyper-path flow to improve network traffic distribution.

At the micro level, the vehicle's route is mapped to match the hyper-path flow. Due to the non-unique mapping between route flows and hyper-path flows, a variety of route flow patterns can be adapted to match the hyper-path flow pattern. Among these potential patterns, we prioritize selecting the pattern that most effectively ensures traffic homogeneity within a region. Consequently, in this approach, users may opt for a sub-optimal choice to meet the above objectives. However, since the routes considered include the current and instantaneously optimal routes, the options available to vehicles are generally favorable.

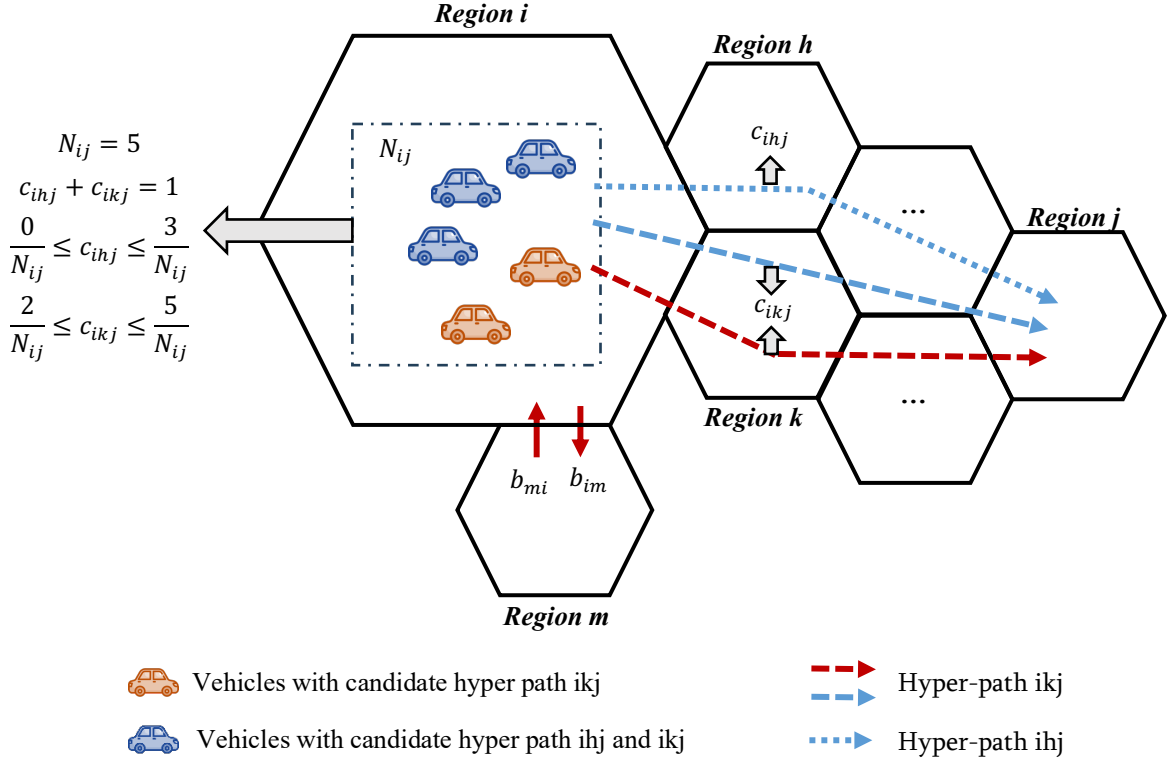

**Fig. 4.** Multi-region network.

III. INTER-REGION CONTROL

*A. Multi-region Network Model*

We utilize an MFD-based model for a multi-region network, incorporating perimeter control and regional route guidance. The discrete-time traffic dynamics for the network are as follows [30, 31]:

$$N_i(t) = N_{ii}(t) + \sum_{j \in R, i \neq j} N_{ij}(t), \forall i \in R \qquad (1)$$

$$N_{ii}(t+1) = N_{ii}(t) + Q_{ii}(t) + \sum_{h \in R_i} N_{hii}(t) - N_i^{II}(t),$$

$$\forall i \in R \qquad (2)$$

$$N_{ij}(t+1) = N_{ij}(t) + Q_{ij}(t) + \sum_{h \in R_i, h \neq j} N_{hij}(t) - \sum_{h \in R_i} N_{ihj}(t),$$

$$\forall i, j \in R, i \neq j \qquad (3)$$

The number of *Type I and Type II vehicles*, which together constitute the completion flow, can be estimated based on the current number of vehicles within the region at macro control step $t$, using the MFD function, as follows:

$$N_{ij}^I(t) = \frac{N_{ij}(t)}{N_i(t)} G_i\big(N_i(t)\big) T_{\mathrm{mac}}, \forall i, j \in R, i \neq j \qquad (4)$$

$$N_i^{II}(t) = \frac{N_{ii}(t)}{N_i(t)} G_i\big(N_i(t)\big) T_{\mathrm{mac}}, \forall i \in R \qquad (5)$$

In this context, the MFD function $G_i\big(N_i(t)\big)$ is modeled as a third-order polynomial to describe the relationship between completion flow and vehicle accumulation[1]:

$$G_i\big(N_i(t)\big) = \beta_{3,i} N_i^3(t) + \beta_{2,i} N_i^2(t) + \beta_{1,i} N_i(t), \forall i \in R \qquad (6)$$

where $\beta_{1,i}, \beta_{2,i}, \beta_{3,i}$ are the estimated parameters.

We use $N_{ihj}(t)$ to represent vehicle movement among regions. *Type I vehicles*, denoted as $N_{ij}^I(t)$, contain all the vehicles that may be included in $N_{ihj}(t)$, but whether they belong to $N_{ihj}(t)$ depends on the applied control measures. Specifically, *Type I vehicles* in region $i$ with destination region $j$ that choose region $h$ as next region to travel and cross the boundary between regions $i$ and $h$ at macro control step $t$ are counted in $N_{ihj}(t)$.

$$N_{ihj}(t) = b_{ih}(t)c_{ihj}(t)N_{ij}^I(t), \forall i, j \in R, i \neq j, h \in R_i \quad (7)$$

Denote $M_{ihj}(t)$ as the flow rate of vehicles in region $i$ with destination region $j$ which cross the boundary between region $i$ and region $h$ at macro control step $t$.

$$M_{ihj}(t) = \frac{N_{ihj}(t)}{T_{\text{mac}}}, \forall i, j \in R, i \neq j, h \in R_i \quad (8)$$

The flow rate between regional boundary can be calculated as

$$M_{ih}(t) = \sum_{j \in R, j \neq i} M_{ihj}(t), \ \forall i \in R, h \in R_i \quad (9)$$

### B. Joint Perimeter Control and Route Guidance Model

To ensure that the trip completion flow rate reaches its maximum, we make the accumulation of vehicles in each region close to its critical value. We let $b_{ih}(t)$ and $c_{ihj}(t)$ be the decision variable the decision variables of joint control. The joint control model is represented as follows.

$$\min_{\substack{b_{ih}(t) \in \boldsymbol{b}(t), \\ c_{ihj}(t) \in \boldsymbol{c}(t)}} \max_i \left(N_i(t+1) - N_i^{\text{crit}}\right) \quad (10)$$

Subject to:

$$\text{Eq. (1)-(7)}$$

$$M_{ih}^{\min}(t) \leq \sum_{j \in R, j \neq i} b_{ih}(t) c_{ihj}(t) \frac{N_{ij}(t)}{N_i(t)} G_i(N_i(t))$$

$$\leq M_{ih}^{\max}(t) f_{ih}(t), \ \forall i \in R, h \in R_i \quad (11)$$

$$\sum_{h \in R_i} c_{ihj}(t) = 1, \forall i, j \in R, i \neq j \quad (12)$$

$$c_{ihj}^{\min}(t) \leq c_{ihj}(t) \leq c_{ihj}^{\max}(t), \forall i, j \in R, i \neq j, h \in R_i \quad (13)$$

Where $\boldsymbol{b}(t) = \{b_{ih}(t) | \forall i \in R, h \in R_i\}$, $\boldsymbol{c}(t) = \{c_{ihj}(t) | \forall i, j \in R, i \neq j, h \in R_i\}$.

Constraint (11) impose minimum and maximum bounds, denoted as $M_{ih}^{\min}(t), M_{ih}^{\max}(t)$, on the number of vehicles passing between any two adjacent regions. These bounds are determined by several factors, including queue lengths, downstream link capacity, and the saturation flow rate. The specific calculation of these values will be provided in Section IV.

Constraint (12) ensures that the sum of the proportion of route choices for each OD equals to 1. Constraint (13) set upper and lower bounds on the proportion of vehicles choosing hyper-paths. Specifically, the maximum value $c_{ihj}^{\max}(t)$ represents the proportion of vehicles in $N_{ij}(t)$ whose possible next regions contain $h$ select $h$ as the next region, while the minimum value $c_{ihj}^{\min}(t)$ is the proportion of vehicles in $N_{ij}(t)$ that are required to take $h$ as the next region. Fig.4 provides an example of how to determine these upper and lower bounds for $c_{ihj}(t)$ based on the number of vehicles in region $i$ to destination region $j$ and their respective candidate hyper-paths.

The original objective function (Eq. 10) can be reformulated as a linear objective function subject to a set of constraints, as shown in Eqs. (14)-(15). In this formulation, an auxiliary variable $z$ is introduced to represent the maximum deviation of the regional accumulations from their critical values.

$$\min_{\substack{b_{ih}(t) \in \boldsymbol{b}(t), \\ c_{ihj}(t) \in \boldsymbol{c}(t), \\ z}} z \quad (14)$$

$$z \geq N_i(t+1) - N_i^{\text{crit}}, \forall i \in R \quad (15)$$

The joint control problem is formulated as a nonlinear programming model with a linear objective function with quadratic constraints, which can be solved using numerical solvers (e.g., scipy.optimize.minimize). This solution yields the optimal values of the parameter $z$, perimeter control variables $b_{ih}(t)$ and the route guidance variables $c_{ihj}(t)$.

## IV. Intra-region Control

The proposed joint control model solves the perimeter control variables and route guidance variables in the network (macro)

level. They cannot be used directly to control and route traffic within a region. This section presents a method to map them to each intersection's signal timing and each vehicle's route choice within a region.

*A. Regional Boundary Control*

The regional boundary control is implemented for the boundary of each pair of adjacent regions (that is, $i,h, \forall i \in R, h \in R_i$). In each micro control step, a multi-phase signal timing plan is executed. To determine this multi-phase plan, we need to seek out the multi-phase plan that meets the joint control requirements by comparing the expected flow rate and the estimated flow rates for each multi-phase plan. Among these multi-phase plans the multi-phase plan with the maximum traffic pressure is activated.

1) **Flow Rate Calculation**

Fig. 5 shows a method for updating the expected flow rate at micro control step $k$, $m_{ih}(k)$, based on the observed actual flow rates from previous steps $\tilde{m}_{ih}(z)$ ($z = 1,2, ..., k-1$). The goal of these adjustments is to ensure that the cumulative number of vehicles passing through the boundary between region $i$ and region $h$ matches the expected flow rate at the macro control step, $M_{ih}(t)$. The updates are crucial for maintaining the consistency of boundary flow at the macro and micro control level.

Given the perimeter control variables $b_{ih}(t)$ and the route guidance variables $c_{ihj}(t)$ solved by the model in Section III, the expected traffic flow rate of macro control step $M_{ih}(t)$ can be calculated by (4)-(9). And the expected flow rate of micro control step can be calculated as

$$m_{ih}(k) = \frac{M_{ih}(t)T_{\text{mac}} - \sum_{z=1}^{z=k-1}\tilde{m}_{ih}(z)T_{\text{mic}}}{(u-k+1)T_{\text{mic}}} \quad (16)$$

where $u$ is the number of micro control steps of each macro control step.

The expected traffic flow may be achieved by setting an eligible multi-phase plan for the boundary. The flow rate for a multi-phase plan can be estimated according to the number of arriving vehicles, saturated flow rate, and downstream link capacity.

| **Algorithm**: Perimeter Control Strategy |
|---|
| 1. Calculate the expected traffic flow rate $m_{ih}(k)$ according to Eq. (16). |
| 2. Calculate the traffic flow rate $\hat{m}_{ih}^s(k)$ for each multi-phase plan $s$ for the current micro control step according to Eq. (17). |
| 3. Determine the set $S_{ih}^m(k)$ of multi-phase plan that comply with the expected flow rate requirements according to Eq. (18). |
| 4. Calculate the weight of these multi-phase plans in $S_{ih}^m(k)$ according to Eq. (22). |
| 5. Select the multi-phase plan $s_{ih}^{\text{act}}(k)$ with the maximum weight as described in Eq. (23). |

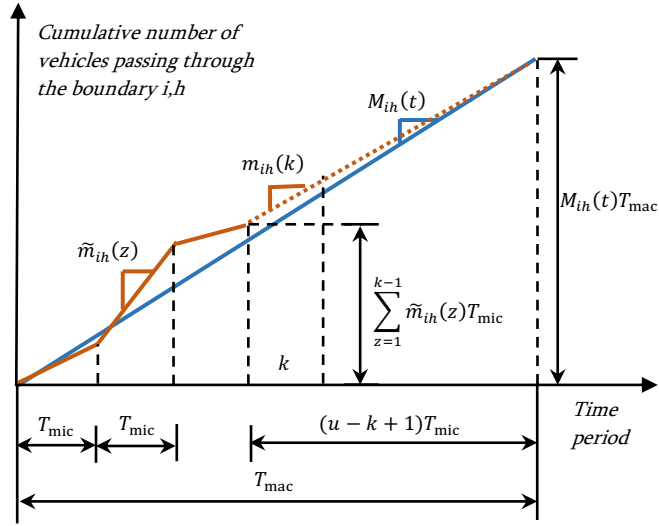

**Fig. 5.** Method for updating expected flow rate at the micro control step.

$$\hat{m}_{ih}^s(k) = \frac{1}{T_{\text{mic}}}\sum_{p\in P_s}\sum_{l\in L_{ih}^p}\min\left\{\begin{matrix}e_l(k), c_lT_{\text{mic}}, \\ \frac{1}{n_l}\sum_{o\in O_l}(l_o - q_o(k))\end{matrix}\right\} \quad (17)$$

where $L_{ih}^p$ is the set of approach lanes from region $i$ to region $h$ which is activated at phase $p$; $O_l$ is the output lane of the input lane $l$; $c_l$ is the saturated flow rate of lane $l$; $e_l(k)$ is the number of vehicles arriving at the stop line or joining the queue in lane $l$ at the time step $k$, which can be obtained based on vehicles' current position and speed; $q_o(k)$ is the queued vehicle number of lane $o$ at time step $k$; $n_l$ is the number of output lanes of the input lane $l$, $l_o$ is the capacity of the lane $o$.

Define $S_{ih}^{\mathrm{m}}(k)$ as the set of multi-phase plan $s$ at micro step $k$ that comply with the expected flow rate requirements at the boundary between region $i$ and region $h$:

$$S_{ih}^{\mathrm{m}}(k) = \left\{ s \left| \begin{array}{c} \frac{|\hat{m}_{ih}^{s}(k) + \hat{m}_{ih}^{\mathrm{ng}}(k) - m_{ih}(k)|}{m_{ih}(k)} < (u - k + 1)\sigma, \\ \frac{|\hat{m}_{hi}^{s}(k) + \hat{m}_{hi}^{\mathrm{ng}}(k) - m_{hi}(k)|}{m_{hi}(k)} < (u - k + 1)\sigma, \\ \forall s \in S_{ih} \end{array} \right. \right\} \quad (18)$$

The estimated flow rate through the non-gating intersection $\hat{m}_{ih}^{\mathrm{ng}}(k)$, can be estimated according to the flow rate of the previous time step. Specifically, $\hat{m}_{ih}^{\mathrm{ng}}(k) = \tilde{m}_{ih}^{\mathrm{ng}}(k-1)$, $\tilde{m}_{ih}^{\mathrm{ng}}(k-1)$ is the observed flow rate from region $i$ to region $h$ through the non-gating intersections at micro control step $k-1$; $\sigma$ is a base threshold of the difference between the expected flow rate and the estimated flow rate under multi-phase plan $s$. $(u - k + 1)\sigma$ is the tolerance threshold at step $k$, which decreases as the step progresses. This design allows the system to begin with boundary flow rates slightly above or below the expected value in earlier steps, and to compensate with opposite deviations in later steps. Such flexibility expands the feasible multi-phase plans and enhances the efficiency and adaptability of boundary flow control.

In addition, the maximum and minimum flow rate from region $i$ to region $h$ at macro control step $t$ can be determined according to the estimated flow rate at the beginning of the macro control step, the first micro control step.

$$M_{ih}^{\max}(t) = \max_{s \in S_{ih}} \hat{m}_{ih}^{s}(0) + \hat{m}_{ih}^{\mathrm{ng}}(0) \quad (19)$$

$$M_{ih}^{\min}(t) = \min_{s \in S_{ih}} \hat{m}_{ih}^{s}(0) + \hat{m}_{ih}^{\mathrm{ng}}(0) \quad (20)$$

### 2) Multi-phase Plan Selection based on Backpressure Policy

Since there are usually several multi-phase plans that satisfy the requirements, we must choose one to activate. A backpressure policy is adopted to select a multi-phase plan considering all the movements at the boundary intersections, not just the movements between two regions.

The multi-phase plan with the maximum traffic pressure is selected at each micro control time step. According to the backpressure control policy, the weight of each phase of each intersection can be calculated as[36]:

$$w_p(k) = \sum_{l \in L_p} \left( q_l(k) - \frac{1}{n_l} \sum_{o \in O_l} q_o(k) \right) c_l \quad (21)$$

where $q_l(k)$ is the queued vehicle number of lane $l$ at time step $k$.

The weight of a multi-phase plan s is as follows,

$$W_s(k) = \sum_{p \in P_s} w_p(k), \forall s \in S_{ih}^{\mathrm{m}}(k) \quad (22)$$

From the set of multi-phase plans $s$ belonging to $S_{ih}^{\mathrm{m}}$, the multi-plan $s_{ih}^{\mathrm{act}}(k)$ that possesses the highest weight is selected for activation for boundary between region $i$ and region $h$ at the micro control step $k$.

$$s_{ih}^{\mathrm{act}}(k) = \mathrm{argmax}_{s \in S_{ih}^{\mathrm{m}}(k)} W_s(k) \quad (23)$$

### 3) Perimeter Control Strategy

The specific steps of the perimeter control are shown in **Algorithm**.

Compared with previous perimeter control, the proposed method adopted a more accurate approach to estimate the flow rate through the boundary of adjacent regions and apply it to the multi-phase plan selection and the constraints of the joint control model. Moreover, movements other than the movements on the regional (macro) hyper-path at the boundary intersections are also considered.

### B. Route choice model

To optimize the route choice of vehicles within a region and achieve balanced traffic flow, the route choice model is defined:

$$\min_{\phi_{ijr}^{v}(k)} \beta \sum_{j \in R, i \neq j} \sum_{h \in R_i} \left( \frac{\sum_{v \in V_{ij}} \sum_{r \in U_v, G_r = h} \phi_{ijr}^{v}(k)}{N_{ij}(k)} - c_{ihj}(t) \right)^2 + \sum_{x \in X_i} (d_x(k) - \bar{d}_i(k))^2 \quad (24)$$

Subject to:

$$\sum_{r \in U_v} \phi_{ijr}^{v}(k) = 1, \forall i, j \in R, i \neq j, v \in V_{ij} \quad (25)$$

$$0 \leq \phi_{ijr}^{v}(k) \leq 1, \forall i, j \in R, i \neq j, v \in V_{ij}, r \in U_v \quad (26)$$

where $k \in K_t$; $\beta$ is a parameter.

The objective function consists of two main components. The first term aims to minimize the difference between the proportion of vehicles that select a route and the solution of joint control model. Specifically, the proportion of vehicles in region $i$ to destination region $j$ with region $h$ as next region to travel can be represented as $\frac{\sum_{v \in V_{ij}} \sum_{r \in U_v, G_r = h} \phi^v_{ijr}(k)}{N_{ij}(k)}$, and it closely matches the route guidance variable $c_{ihj}(t)$.

The second term $\sum_{x \in X_i}(d_x(k) - \bar{d}_i(k))^2$ ensures the homogeneity of traffic conditions inside the region, implying that the estimated link density for each link in region $i$, noted as $d_x(k)$, should be close to average link density of region $i$, $\bar{d}_i(k)$. The expressions for $d_x(k)$ and $\bar{d}_i(k)$ are given by the following equations:

$$d_x(k) = \frac{\sum_{j \in R} \sum_{v \in V_{ij}} \sum_{L_r = x} \phi^v_{ijr}(k)}{l^a_x l^e_x}, \forall x \in X_i, i \in R \qquad (27)$$

$$\bar{d}_i(k) = \frac{N_i(k)}{\sum_{x \in X_i} l^a_x l^e_x}, \forall i \in R \qquad (28)$$

where $l^a_x$, $l^e_x$ are the lane number and link length of link $x$, and $L_r$ is the link where the vehicle is estimated to be at the end of the $k$th time step. It is assumed that vehicles coming from or going to other areas during the current control step are ignored when calculating the link density.

Moreover, constraint (25) ensures that the sum of the route proportions is equal to 1 for each vehicle, and constraint (26) restricts the route proportions to be between 0 and 1.

The route choice model is also a convex quadratic programming model characterized by a convex quadratic objective function and linear constraints.

## V. Case Study

### A. Network Description

To evaluate the performance of the proposed control framework, we take a road network of Yangzhou for our case study, as shown in Fig. 6. The network consists of 1102 links (one-way) and 347 intersections, 149 of which are signalized.

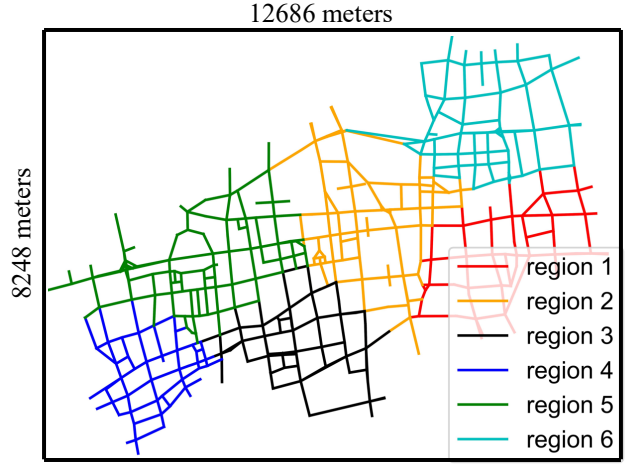

**Fig. 6.** Simulated road network

Simulation of Urban Mobility (SUMO) is employed to simulate traffic conditions in our case study. An application tool of SUMO, jtcrouter.py is used to generate traffic demand that produces congestion in the network for a period of time during the peak hour.

The source codes and simulation configuration files are available at: https://github.com/xypeng12/Multi-Region-Traffic-Control-with-SUMO.

Since the MFD model relies on homogeneity within each region to ensure a unimodal, low-scatter, and demand-insensitive relationship between vehicle accumulation and trip completion rates [2, 37], we applied the RatioCut algorithm to partition the network into 6 homogeneous regions [38]. This method minimizes the variance in link density, allowing us to obtain well-fitted MFD functions for each region.

To derive the MFD for each region, we collected accumulation and outflow data from simulations lasting 7200 seconds under varying demand levels (25%, 50%, 75%, 100%, and 125% of the base demand). Trip completion flow was approximated by the observed regional outflow. To mitigate the effects of short-term signal fluctuations, data were sampled every 120 seconds. Based on these observations, smooth MFD curves were fitted for each region to capture the relationship between accumulation and trip completion flow.

Table II summarizes the MFD functions and their critical accumulation values. Fig. 7 shows the MFD scatter plots with the fitted curves for each region. To further analyze the spatial heterogeneity and causes of MFD scatter, we computed the coefficient of variation (CV) of density across all roads within each region. As illustrated by the color gradient in Fig. 7, a substantial number of data points with high CV values tend to lie below the fitted curve, indicating that heterogeneity reduces MFD accuracy. Specifically, it may lead to an overestimation of trip completion flow while having little impact on the critical accumulation value.

Nevertheless, this deviation has limited impact on control performance for the following reasons. Although such overestimation introduces a discrepancy between the actual and modeled $G_i\big(N_i(t)\big)$, the boundary flow $M_{ih}(t)$ in our control model is further constrained by local traffic information through constraint (11). These upper and lower bounds prevent $M_{ih}(t)$ from reaching overestimated levels, even when $G_i\big(N_i(t)\big)$ is inaccurate. This design helps mitigate the impact of MFD deviations and enhances control robustness. Moreover, since our controller operates under an MPC framework and continuously updates decisions based on real-time traffic conditions, it remains effective even with imperfect MFD estimation.

### B. Joint Control Settings

#### 1) Scenarios

To illustrate the effect of multi-region joint control and to examine the roles of vehicle routing and boundary signal control, we compare several control strategies built upon our proposed multi-region framework, along with baseline strategies such as boundary backpressure control. These strategies differ in whether and how they incorporate perimeter control and route guidance.

• The proposed multi-scale joint control of perimeter

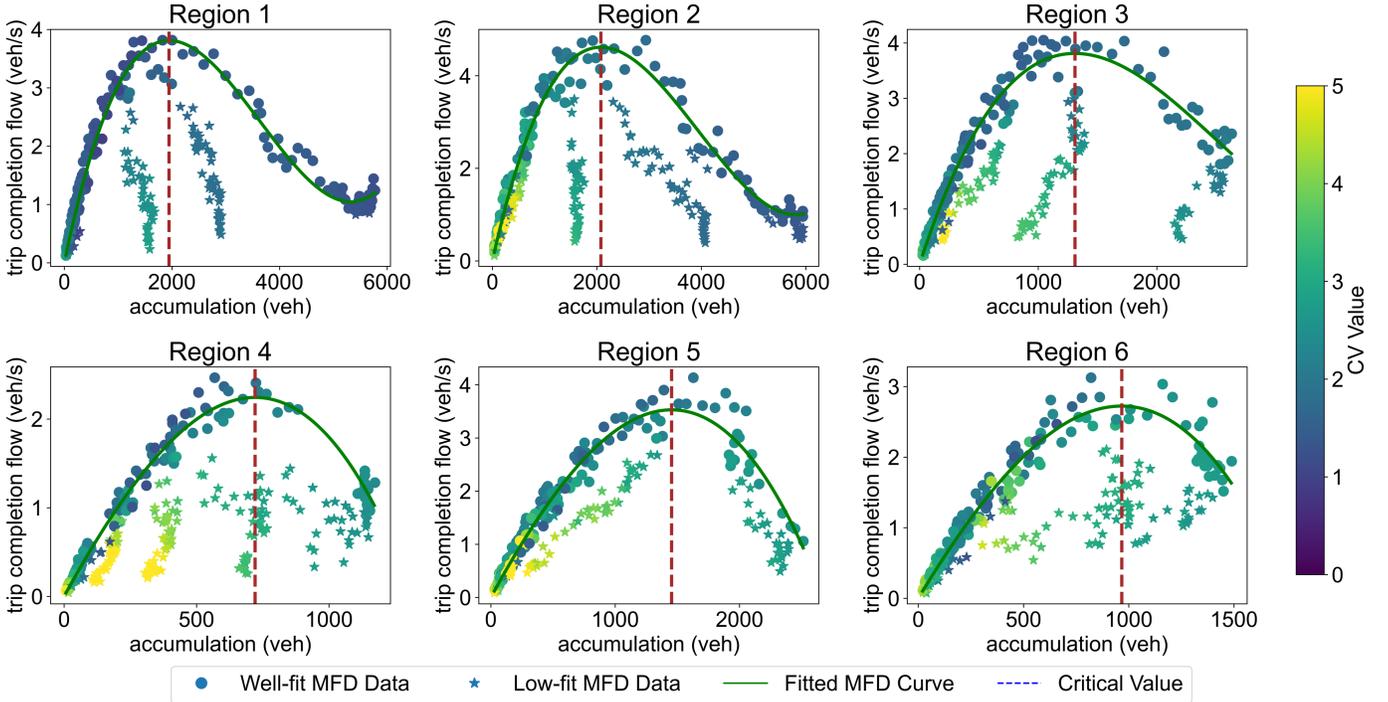

**Fig. 7.** MFDs of each region.

TABLE II
Mfd Function and Critical Accumulation

| Region | MFD function | Critical accumultion (veh) |
|---|---|---|
| 1 | $G_1(N_1(t)) = 1.44 \times 10^{-10} N_1^3(t) - 1.57 \times 10^{-6} N_1^2(t) + 4.46 \times 10^{-3} N_1(t)$ | 1946 |
| 2 | $G_2(N_2(t)) = 1.39 \times 10^{-10} N_2^3(t) - 1.65 \times 10^{-6} N_2^2(t) + 5.04 \times 10^{-3} N_2(t)$ | 2077 |
| 3 | $G_3(N_3(t)) = 4.50 \times 10^{-10} N_3^3(t) - 3.40 \times 10^{-6} N_3^2(t) + 6.59 \times 10^{-3} N_3(t)$ | 1310 |
| 4 | $G_4(N_4(t)) = -1.46 \times 10^{9} N_4^3(t) - 2.21 \times 10^{-6} N_4^2(t) + 5.46 \times 10^{-3} N_4(t)$ | 721 |
| 5 | $G_5(N_5(t)) = -2.59 \times 10^{-10} N_5^3(t) - 9.18 \times 10^{-7} N_5^2(t) + 4.31 \times 10^{-3} N_5(t)$ | 1454 |
| 6 | $G_6(N_6(t)) = -7.38 \times 10^{-10} N_6^3(t) - 1.49 \times 10^{-6} N_6^2(t) + 4.95 \times 10^{-3} N_6(t)$ | 967 |

control and route guidance, MSJC: This is the fully integrated strategy proposed in this paper. Given regional-level traffic states $\{N_{ij}(t), Q_{ij}(t)\}$ obtained from SUMO, the multi-region model is constructed to determine the upper-level control variables $\{M_{ih}(t), c_{ihj}(t)\}$. At the lower level, SUMO simulates individual vehicle movements and enables the control of each vehicle's routes $\phi_{ijr}^v(k)$ and boundary intersection signals $s_{ih}^{\text{act}}(k)$ according to the lower-level control strategies.

• The multi-scale perimeter control with a logit-based route guidance strategy, MSPC-LR: this strategy is inspired by the proportional-integral(PI)-based perimeter control [3, 18] with logit-based route guidance as proposed in [33]. However, such control logic cannot be directly implemented in a microscopic simulation without additional mechanisms. To ensure applicability and fair comparison, two key components from our proposed framework are integrated. First, the target flow $M_{ih}(t)$ computed by the PI controller is realized through our lower-level control method described in Section IV.A. Second, instead of using a predefined path set as in [33], the route set is dynamically generated via the column generation approach introduced in Section II.B, with consistent routing flexibility.

• The backpressure control for boundary intersections with the logit-based route guidance strategy, BP-LR: This strategy applies the backpressure control method [36] at boundary intersections and incorporates logit-based route guidance. As in MSPC-LR, logit-based route choice is supported by our column generation procedure

• The multi-scale perimeter control without route guidance, MSPC: This configuration only includes perimeter control. The upper-level flow $M_{ih}(t)$ is determined using the PI controller and its implementation relies on the lower-level control method described in Section IV.A. No route adjustment is performed.

• The backpressure control for boundary intersections without route guidance, BP: this baseline applies adaptive backpressure control [36] at boundary intersection without any routing module.

Among these strategies, MSJC represents the complete proposed method, while MSPC-LR, MSPC and BP-LR are partially integrated variants designed for applicability and fair comparison. BP serves as baseline configuration for comparison.

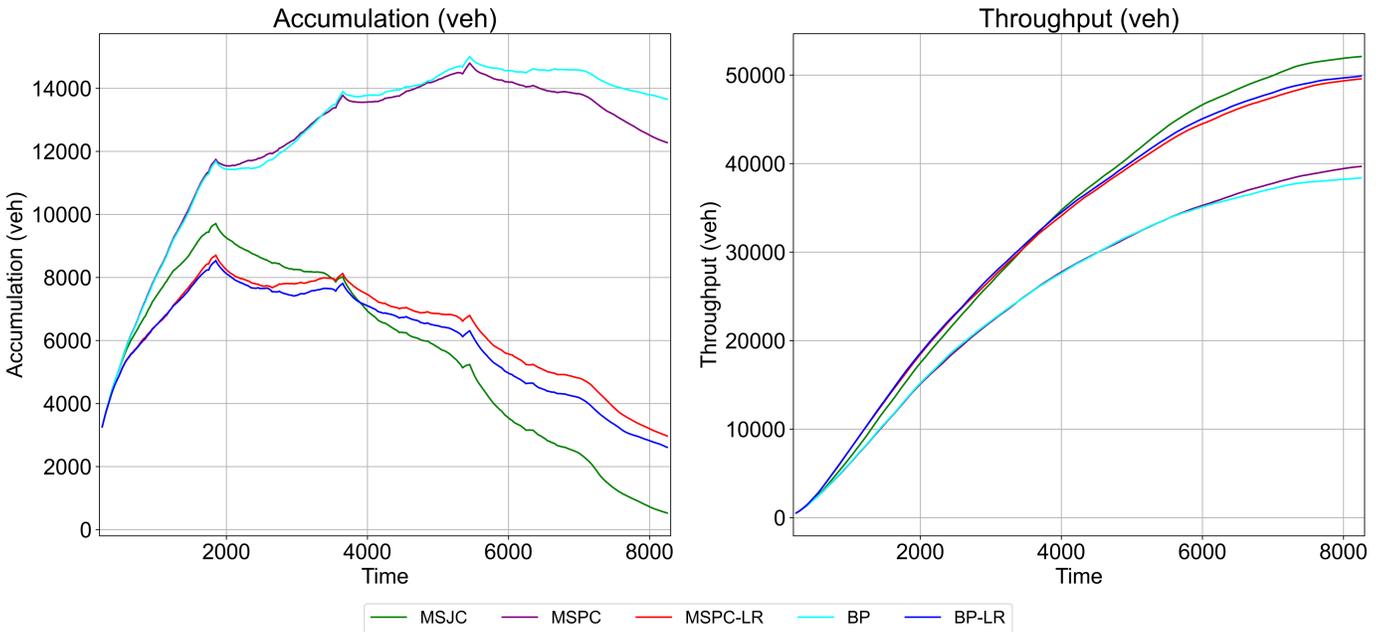

**Fig. 8.** Network accumulation and accumulated throughput.

TABLE III
OVERALL PERFORMANCE OF FIVE CONTROL STRATEGIES

|  | MSJC | BP-LR | MSPC-LR | MSPC | BP |
|---|---|---|---|---|---|
| Total travel time (in $10^7$ veh·s) | 4.536 | 4.829 | 5.099 | 9.914 | 10.130 |
| Throughput (in $10^4$ veh) | 5.212 | 4.992 | 4.961 | 3.972 | 3.842 |

### 2) Simulation and control setup

The traffic demand consists of vehicles entering the network in a time window of [0, 7200] seconds. The macro control step length is set to 100 seconds, and the micro control step length is 10 seconds. A 200-second warm-up period is applied at the beginning of the simulation. To observe the full process of traffic growth and dissipation, we simulate the traffic until the network is cleared. Control actions are activated once the accumulation of any region exceeds 30% of its critical value.

Other parameters required for control (e.g., the layout of the network and boundary intersections, current accumulation, queue length, and potential routes for vehicles) are obtained from the simulation.

### C. Results and Analysis

Fig. 8 represents the overall accumulation and accumulated throughput over time, reflecting the performance of the five control strategies. As summarized in Table III, the strategies are ranked based on both total travel time and throughput as follows: MSJC, BP-LR, MSPC-LR, MSPC, and BP. This trend is also reflected in the accumulation and cumulative throughput, as shown in Fig.9 and Fig. 10. Fig.11 illustrates the MFD scatter plots for multiple regions. Under MSJC control, most regions (except Region 3) exhibit MFD points that closely follow the measured MFD curve, supporting the validity of the assumption underlying the proposed method. Similar patterns can be observed under MSPC-LR and BP-LR; however, the points for Region 1 and Region 2 fall below the MFD curve due to heterogeneous traffic conditions—specifically, localized gridlock within these regions.

Under MSJC, vehicle accumulation in each region approaches the critical value, thereby optimizing throughput and demonstrating significant improvement over other control methods. Specifically, during the congested period from approximately 1600 to 4000 seconds, the accumulations in Regions 1, 2, 5, and 6 are consistently maintained near their critical values. In contrast, Region 4, which has relatively low demand, sees a higher accumulation under MSJC than under MSPC-LR and BP-LR, as vehicles are encouraged to pass through this region to relieve congestion elsewhere. Region 3 slightly exceeds its critical accumulation, possibly due to highly heterogeneous traffic conditions and an overestimation of the MFD as shown in Fig.11. However, despite the excess accumulation, its throughput remains comparable to that under MSPC-LR and BP-LR, indicating that its performance is still acceptable.

MSPC-LR also performs well across multiple regions but fails to control inflows into Region 1, leading to some congestion. BP-LR performs similarly to MSPC-LR.

Routing in MSJC is a first-best route choice aimed at improving network performance, while the real-time logit-based route guidance strategies in MSPC-LR and BP-LR offer second-best routing. Though these methods guide vehicles to less congested areas, their vehicle-level perspective limits overall effectiveness compared to network-wide strategies.

MSPC does not perform as well as MSJC or MSPC-LR, as it does not guide vehicles toward less congested areas. BP is the worst-performing strategy, as it only adjusts signal timings at boundaries, which is not capable of managing network distribution effectively. Consequently, vehicle accumulation in each region exceeds optimal levels, leading to more congestion and lower throughput. These results highlight the efficiency of the proposed MSJC method and the importance of jointly optimizing signal control and route guidance.

To examine the integration of perimeter control and route guidance, we captured the traffic distribution at the 3200th second and the joint control solutions from 3200 to 3300 seconds, as shown in Fig.12. The number in each region indicates its current accumulation $N_i(3200)$ and the critical accumulation $N_i^{\text{crit}}$. Region 1 and Region 3 are oversaturated, as their current accumulation exceed the critical value, while Region 6 is near-saturated and Region2, Region 4 and Region 5 are under-saturated.

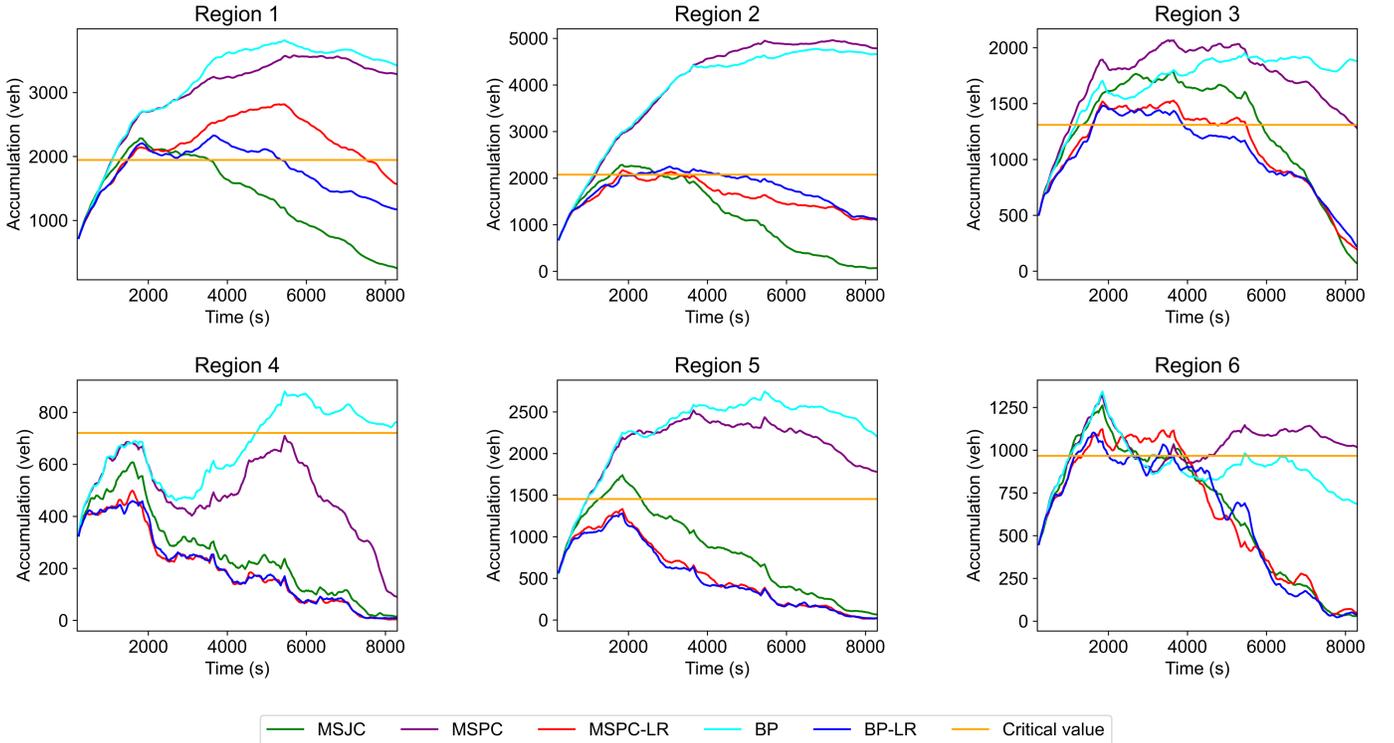

**Fig. 9.** Accumulation dynamics for multiple regions.

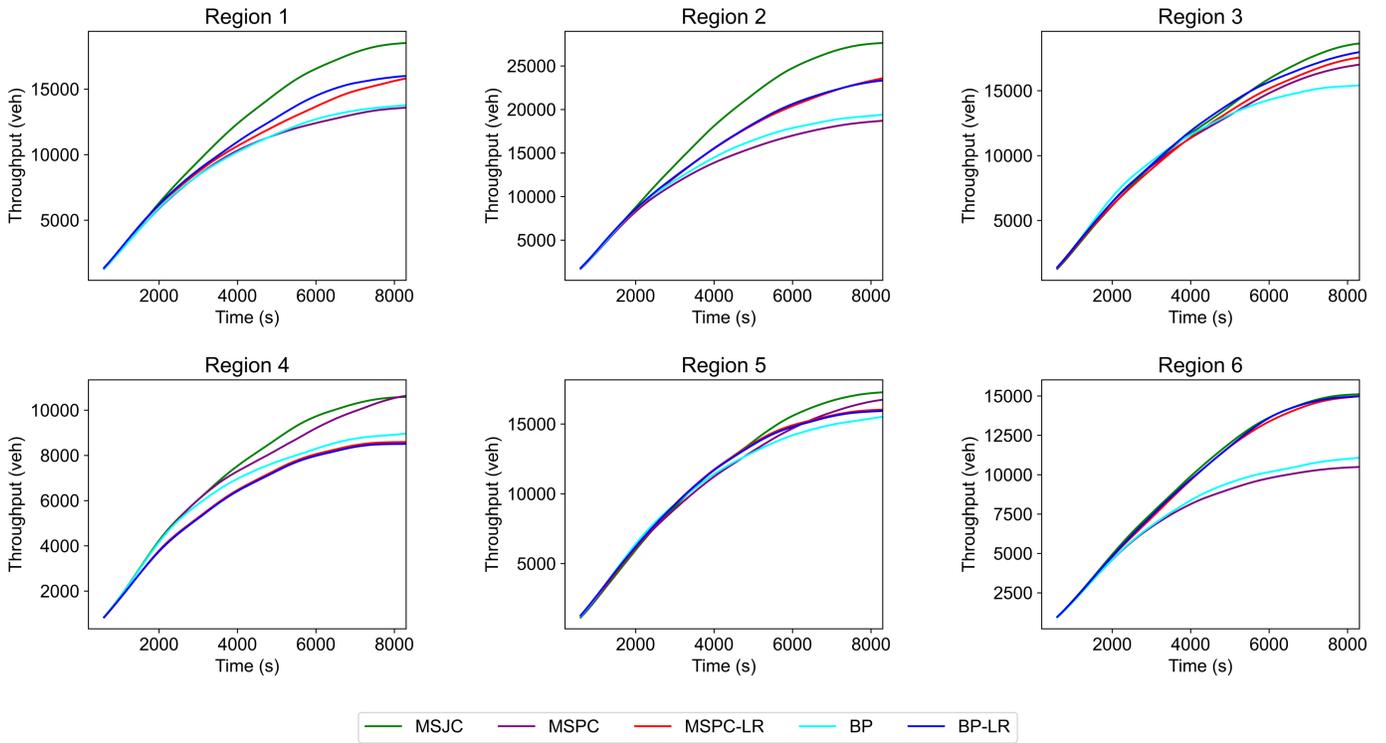

**Fig. 10.** Cumulative throughput for multiple regions.

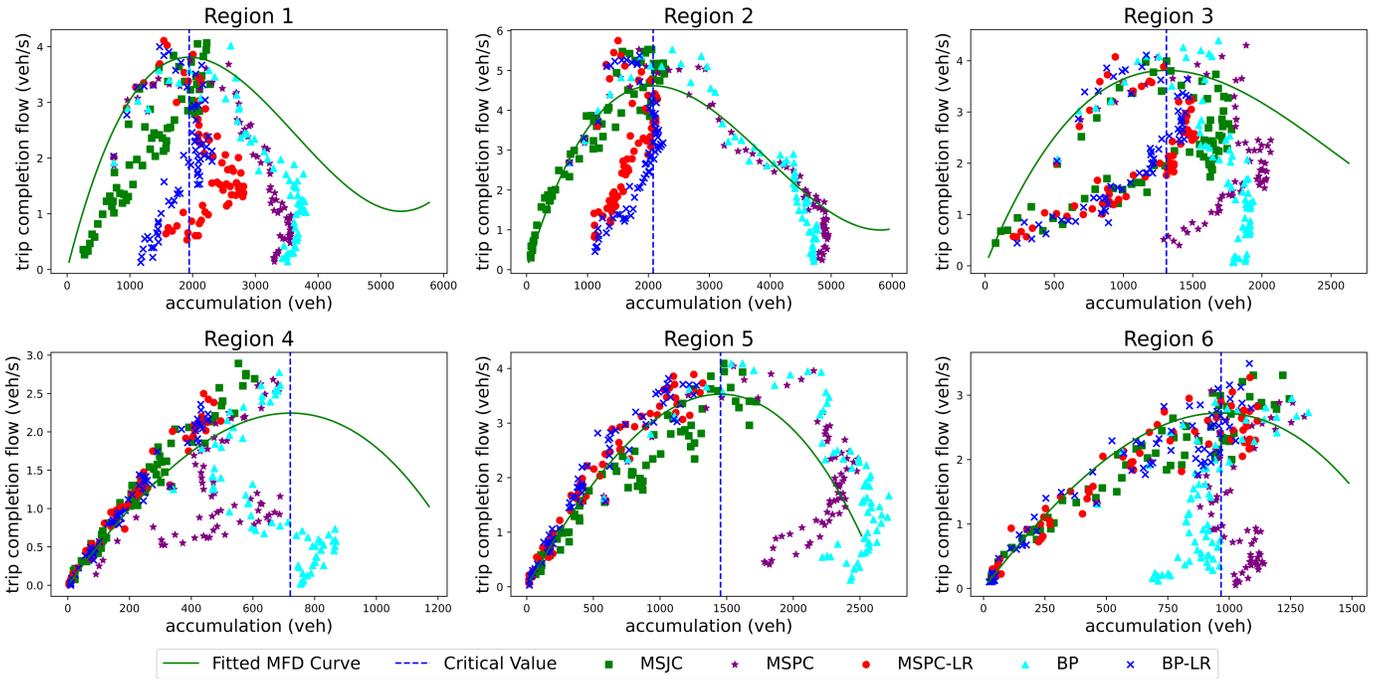

**Fig. 11.** MFDs of multiple regions under control strategies.

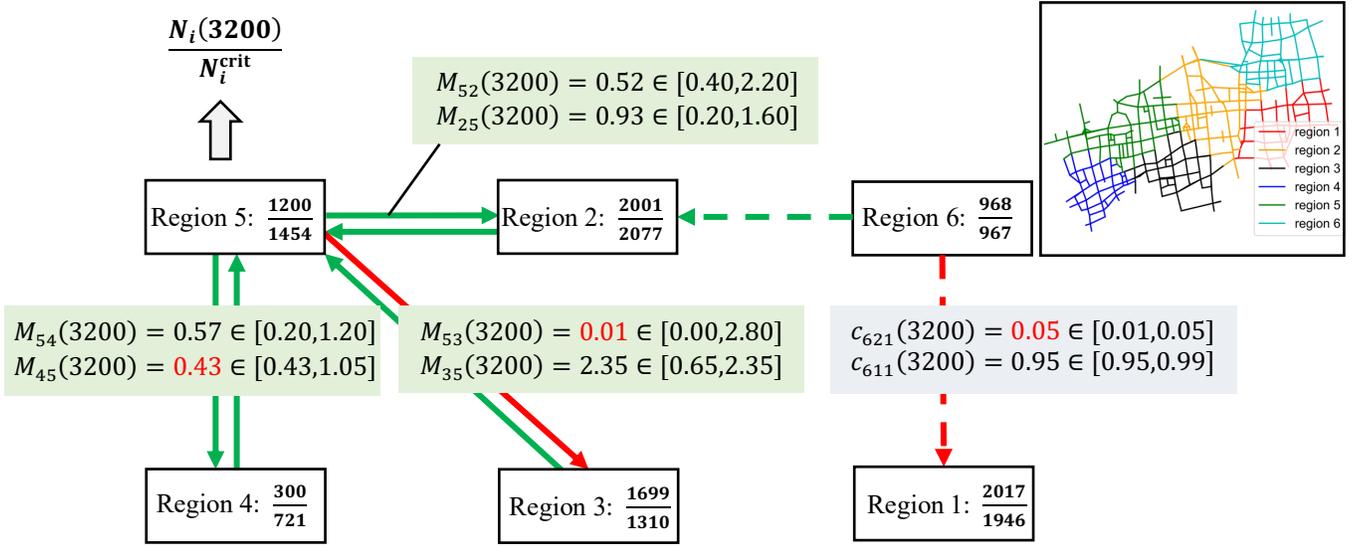

**Fig. 12.** The traffic distribution (3200s) and the control variables (3200s-3300s).

The traffic movement between Region 5 and its adjacent regions is depicted, where flow from lower-density regions to higher-density regions (e.g. $M_{53}(3200)$, $M_{45}(3200)$) is restricted, with the flow rate approaching its minimum value. Conversely, traffic flow from oversaturated region to undersaturated region, that is, $M_{35}(3200)$ reaches its maximum. Meanwhile, vehicles are guided towards less congested areas. For instance, vehicles in Region 6 are directed to the undersaturated Region 2 as much as possible. However, the magnitude of this flow is limited because most vehicles' available routes from Region 6 lead directly through the oversaturated Region 1. In summary, the proposed joint control method effectively manages the inflow into regions with high traffic loads and redistributes traffic across the network. However, its effectiveness is limited by regional boundary conditions, such as the number of approaching vehicles, the saturated flow rate, and the downstream capacity. The flexibility in guiding vehicles depends on the vehicles within the region and their available candidate routes.

To validate that the detailed control actions align with the joint control solution, we compared the expected flow rates solved by the macro-level control model with the actual flow rates observed from the simulation under the specific micro-level control measures. The boundary between Region 2 and Region 5 serves as an example, as shown in Fig. 13. The flow from Region 5 to Region 2 consistently matches the expected values throughout the simulation. The flow from Region 2 to Region 5 also aligns well most of the time. However, during the period from 2500 to 3100 seconds, a noticeable mismatch occurs. This discrepancy arises because the joint control strategy at that time aims to restrict inflow from the lower-density Region 5 to the higher-density Region 2. Since traffic signals typically operate with symmetric green phases, it is difficult to significantly limit inflow in one direction while maintaining high outflow in the other. In other words, few multiphase signal plans can meet the strict thresholds set by the upper-level control objectives. Overall, the expected and actual flow rates showed similar patterns. The observed differences are acceptable due to necessary compromises in both directions.

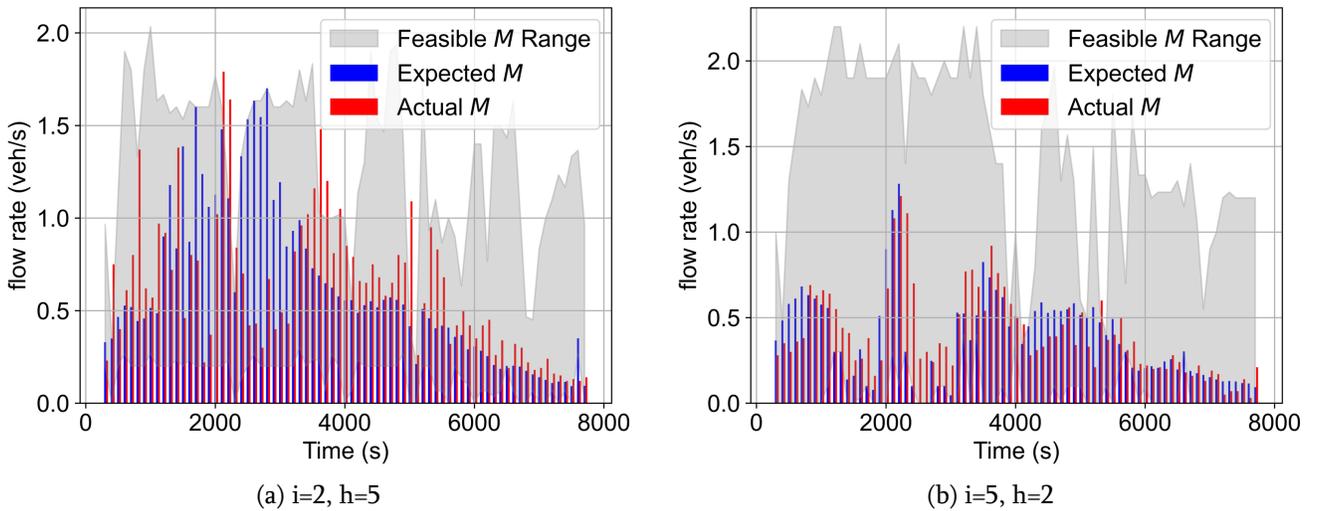

(a) i=2, h=5

(b) i=5, h=2

**Fig. 13.** Flow rates between two adjacent regions.

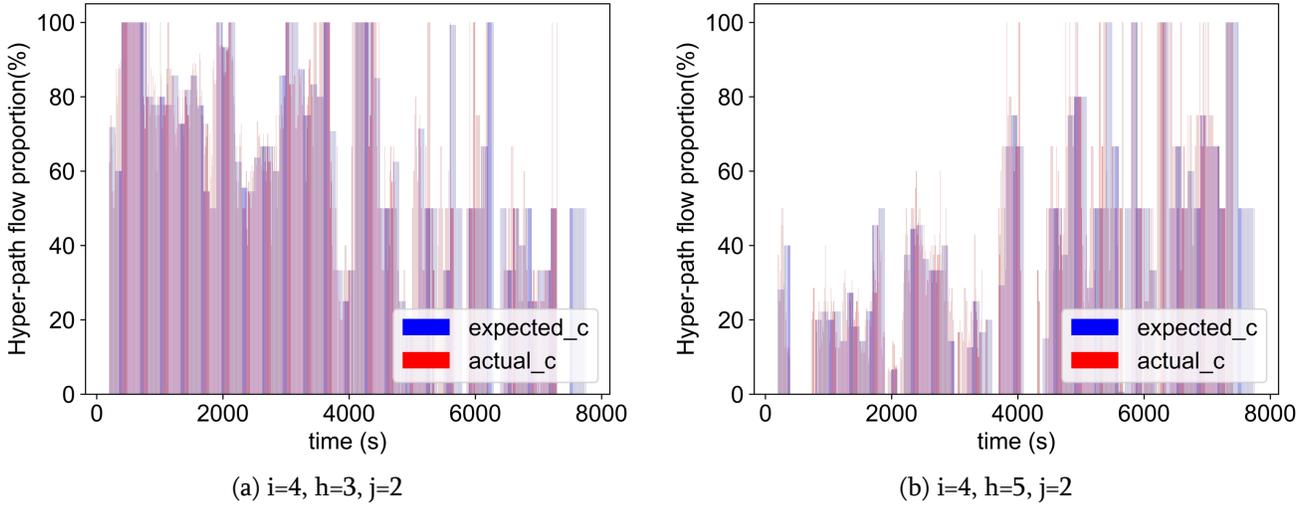

(a) i=4, h=3, j=2
(b) i=4, h=5, j=2

**Fig. 14.** Hyper-path flow proportions.

Fig. 13 also presents the feasible ranges of the boundary flow rates, such as $\left[M_{25}^{\min}(t), M_{25}^{\max}(t)\right]$ and $\left[M_{52}^{\min}(t), M_{52}^{\max}(t)\right]$. Both the expected and actual flow rates consistently fall within these ranges, demonstrating the effectiveness of the boundary flow prediction method.

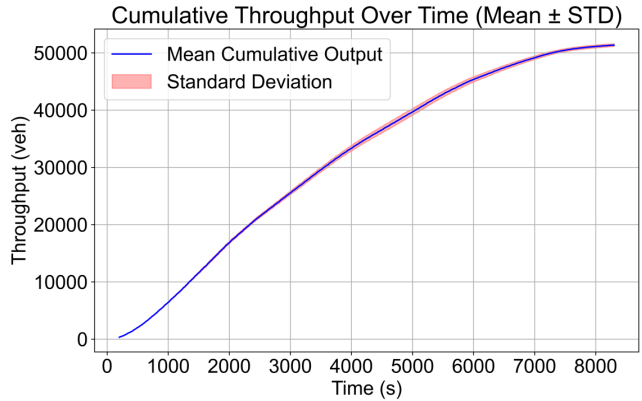

**Fig. 15.** Cumulative throughput over time across 10 replications

Similarly, we examined the alignment between expected and actual proportions of hyper-path flows. Fig. 14 shows the proportions of vehicles in Region 4, with a destination in Region 2, choosing Region 3 or Region 5 as the next region. The expected region-based route proportions were achieved by assigning routes to vehicles, and the actual route flow proportions closely matched the expected ones. After 4800 seconds, when Region 5 became low-density and Region 3 remained oversaturated (as shown in Fig. 9), vehicles were guided to Region 5 as much as possible.

To take into account the stochasticity in SUMO simulations, we conducted 10 replications of the experiment using different random seeds to account for variability. Fig. 15 shows the cumulative throughput over time, where the blue line represents the mean cumulative output, and the red shaded area indicates the standard deviation across the 10 replications. The narrow standard deviation band indicates limited variability in the cumulative throughput, confirming the consistency and reliability of the results.

## VI. Conclusion

This paper proposes a multi-scale traffic control and route guidance system for a large-scale traffic network. At the macroscopic level, a joint control model is proposed to manage the accumulation of each region to ensure the maximum throughput and addresses the inflow rate between adjacent regions and the hyper-path flows for the network. To apply these solutions of the joint model to a specific control scheme for the real world, regional boundary control and route choice model are constructed respectively at the microscopic level. Regional boundary control is carried out through the help of the backpressure algorithm at the gating intersections on the regional boundary. Among several multi-phase plans that meet the requirement of the inflow rates solved by the joint control model, the multi-phase plan with maximum pressure is selected to be activated so that all the movements of the boundary intersections are considered. At the same time, the route choice model solves the actual travel routes for each vehicle to achieve the desired hyper-path flows determined by the joint model and ensures intra-region traffic homogeneity.

The proposed control system is applied and evaluated on a traffic network in Yangzhou, China. The simulation results show that MSJC outperforms other control strategies in controlling regional accumulation and improving throughput.

There are still several issues that deserve further investigation. Excessive intra-region demand may limit the effectiveness of inflow regulation. In such cases, demand-side adjustments (e.g., congestion pricing) could be introduced to further enhance system performance. Another issue that warrants attention is the discrepancy between the expected and realized flows through some regions. This could potentially be addressed by better aligning routes with gating intersections within the region.

## Reference

[1] H. Ding, F. Guo, X. Zheng, and W. Zhang, "Traffic guidance-perimeter control coupled method for the congestion in a macro network," *Transportation Research Part C-Emerging Technologies,* vol. 81, pp. 300-316, Aug 2017, doi: 10.1016/j.trc.2017.06.010.

[2] C. F. Daganzo and N. Geroliminis, "An analytical approximation for the macroscopic fundamental diagram of urban traffic," *Transportation Research Part B: Methodological,* vol. 42, no. 9, pp. 771-781, 2008.

[3] M. Keyvan-Ekbatani, A. Kouvelas, I. Papamichail, and M. Papageorgiou, "Exploiting the fundamental diagram of urban networks for feedback-based gating," *Transportation Research Part B: Methodological,* vol. 46, no. 10, pp. 1393-1403, 2012.

[4] J. Haddad and A. Shraiber, "Robust perimeter control design for an urban region," *Transportation Research Part B: Methodological,* vol. 68, pp. 315-332, 2014.

[5] N. Geroliminis, J. Haddad, and M. Ramezani, "Optimal Perimeter Control for Two Urban Regions With Macroscopic Fundamental Diagrams: A Model Predictive Approach," *IEEE Transactions on Intelligent Transportation Systems,* vol. 14, no. 1, pp. 348-359, 2013, doi: 10.1109/TITS.2012.2216877.

[6] F. de Souza, R. M. Saucedo, O. Mousavizadeh, R. C. Carlson, and M. Keyvan-Ekbatani, "On the evaluation and selection of network-level traffic control policies: Perimeter control, TUC, and their combination," *Transportation Research Part A: Policy and Practice,* vol. 186, p. 104161, 2024.

[7] S. Jiang and M. Keyvan-Ekbatani, "Hybrid perimeter control with real-time partitions in heterogeneous urban networks: An integration of deep learning and MPC," *Transportation Research Part C: Emerging Technologies,* vol. 154, p. 104240, 2023.

[8] J. Haddad and B. Mirkin, "Adaptive perimeter traffic control of urban road networks based on MFD model with time delays," *International Journal of Robust and Nonlinear Control,* vol. 26, no. 6, pp. 1267-1285, 2016, doi: https://doi.org/10.1002/rnc.3502.

[9] K. Ampountolas, N. Zheng, and N. Geroliminis, "Macroscopic modelling and robust control of bi-modal multi-region urban road networks," *Transportation Research Part B: Methodological,* vol. 104, pp. 616-637, 2017/10/01/ 2017, doi: https://doi.org/10.1016/j.trb.2017.05.007.

[10] R. X. Zhong, C. Chen, Y. P. Huang, A. Sumalee, W. H. K. Lam, and D. B. Xu, "Robust perimeter control for two urban regions with macroscopic fundamental diagrams: A control-Lyapunov function approach," *Transportation Research Part B: Methodological,* vol. 117, pp. 687-707, 2018/11/01/ 2018, doi: https://doi.org/10.1016/j.trb.2017.09.008.

[11] J. Wang, X. He, S. Peeta, and X. Yang, "Feedback perimeter control with online estimation of maximum throughput for an incident-affected road network," *Journal of Intelligent Transportation Systems,* vol. 26, no. 1, pp. 81-99, Dec 11 2021, doi: 10.1080/15472450.2020.1797501.

[12] Y. Bichiou, M. Elouni, H. M. Abdelghaffar, and H. A. Rakha, "Sliding Mode Network Perimeter Control," *Ieee Transactions on Intelligent Transportation Systems,* vol. 22, no. 5, pp. 2933-2942, May 2021, doi: 10.1109/tits.2020.2978166.

[13] H. Ding, Y. Di, Z. Feng, W. Zhang, X. Zheng, and T. Yang, "A perimeter control method for a congested urban road network with dynamic and variable ranges br," *Transportation Research Part B-Methodological,* vol. 155, pp. 160-187, Jan 2022, doi: 10.1016/j.trb.2021.11.008.

[14] K. Aboudolas and N. Geroliminis, "Perimeter and boundary flow control in multi-reservoir heterogeneous networks," *Transportation Research Part B: Methodological,* vol. 55, pp. 265-281, 2013/09/01/ 2013, doi: https://doi.org/10.1016/j.trb.2013.07.003.

[15] O. Mousavizadeh and M. Keyvan-Ekbatani, "Regional Transfer Flow Ratio Estimation in a Multi-Reservoir Traffic System," in *2023 IEEE 26th International Conference on Intelligent Transportation Systems (ITSC),* 2023: IEEE, pp. 4108-4113.

[16] M. Keyvan-Ekbatani, R. C. Carlson, V. L. Knoop, and M. Papageorgiou, "Optimizing distribution of metered traffic flow in perimeter control: Queue and delay balancing approaches," *Control Engineering Practice,* vol. 110, p. 104762, 2021.

[17] K. Yang, N. Zheng, and M. Menendez, "Multi-scale perimeter control approach in a connected-vehicle environment," *Transportation research procedia,* vol. 23, pp. 101-120, 2017.

[18] D. Tsitsokas, A. Kouvelas, and N. Geroliminis, "Two-layer adaptive signal control framework for large-scale dynamically-congested networks: Combining efficient Max Pressure with Perimeter Control," *Transportation Research Part C: Emerging Technologies,* vol. 152, p. 104128, 2023.

[19] M. Ramezani, J. Haddad, and N. Geroliminis, "Dynamics of heterogeneity in urban networks: aggregated traffic modeling and hierarchical control," *Transportation Research Part B-Methodological,* vol. 74, pp. 1-19, Apr 2015, doi: 10.1016/j.trb.2014.12.010.

[20] H. Fu, N. Liu, and G. Hu, "Hierarchical perimeter control with guaranteed stability for dynamically coupled heterogeneous urban traffic," *Transportation Research Part C: Emerging Technologies,* vol. 83, pp. 18-38, 2017.

[21] L. Leclercq and N. Geroliminis, "Estimating MFDs in Simple Networks with Route Choice," *Procedia - Social and Behavioral Sciences,* vol. 80, pp. 99-118, 2013/06/07/ 2013, doi: https://doi.org/10.1016/j.sbspro.2013.05.008.

[22] V. V. Gayah and C. F. Daganzo, "Effects of Turning Maneuvers and Route Choice on a Simple Network," *Transportation Research Record,* vol. 2249, no. 1, pp. 15-19, 2011, doi: 10.3141/2249-03.

[23] M. Hajiahmadi, V. L. Knoop, B. De Schutter, H. Hellendoorn, and Ieee, "Optimal Dynamic Route Guidance: A Model Predictive Approach Using the Macroscopic Fundamental Diagram," in *16th International IEEE Conference on Intelligent Transportation Systems (ITSC),* The Hague, NETHERLANDS, 2013

Oct 06-09 2013, in IEEE International Conference on Intelligent Transportation Systems-ITSC, 2013, pp. 1022-1028. [Online]. Available: <Go to ISI>://WOS:000346481000164. [Online]. Available: <Go to ISI>://WOS:000346481000164

[24] C. Menelaou, S. Timotheou, P. Kolios, C. G. Panayiotou, and Ieee, "Joint route guidance and demand management for multi-region traffic networks," in *18th European Control Conference (ECC),* Naples, ITALY, 2019

Jun 25-28 2019, 2019, pp. 2183-2188. [Online]. Available: <Go to ISI>://WOS:000490488302034. [Online]. Available: <Go to ISI>://WOS:000490488302034

[25] C. Menelaou, S. Timotheou, P. Kolios, and C. G. Panayiotou, "Joint Route Guidance and Demand Management for Real-Time Control of Multi-Regional Traffic Networks," *IEEE Transactions on Intelligent Transportation Systems,* vol. 23, no. 7, pp. 8302-8315, 2022, doi: 10.1109/tits.2021.3077870.

[26] M. Yildirimoglu and N. Geroliminis, "Approximating dynamic equilibrium conditions with macroscopic fundamental diagrams," *Transportation Research Part B-Methodological,* vol. 70, pp. 186-200, Dec 2014, doi: 10.1016/j.trb.2014.09.002.

[27] M. Yildirimoglu, M. Ramezani, and N. Geroliminis, "Equilibrium Analysis and Route Guidance in Large-Scale Networks with MFD Dynamics," in *21st International Symposium on Transportation and Traffic Theory (ISTTT)*, Kobe, JAPAN, 2015 Aug 05-07 2015, vol. 9, in Transportation Research Procedia, 2015, pp. 185-204, doi: 10.1016/j.trpro.2015.07.011. [Online]. Available: <Go to ISI>://WOS:000380485800006

[28] M. Yildirimoglu, I. I. Sirmatel, and N. Geroliminis, "Hierarchical control of heterogeneous large-scale urban road networks via path assignment and regional route guidance," *Transportation Research Part B: Methodological,* vol. 118, pp. 106-123, 2018, doi: 10.1016/j.trb.2018.10.007.

[29] S. Jiang, C. Q. Tran, and M. Keyvan-Ekbatani, "Regional route guidance with realistic compliance patterns: Application of deep reinforcement learning and MPC," *Transportation Research Part C: Emerging Technologies,* vol. 158, p. 104440, 2024.

[30] I. I. Sirmatel and N. Geroliminis, "Economic Model Predictive Control of Large-Scale Urban Road Networks via Perimeter Control and Regional Route Guidance," *Ieee Transactions on Intelligent Transportation Systems,* vol. 19, no. 4, pp. 1112-1121, Apr 2018, doi: 10.1109/tits.2017.2716541.

[31] Z. Hou and T. Lei, "Constrained Model Free Adaptive Predictive Perimeter Control and Route Guidance for Multi-Region Urban Traffic Systems," *Ieee Transactions on Intelligent Transportation Systems,* vol. 23, no. 2, pp. 912-924, Feb 2022, doi: 10.1109/tits.2020.3017351.

[32] H. Fu, S. Chen, K. Chen, A. Kouvelas, and N. Geroliminis, "Perimeter Control and Route Guidance of Multi-Region MFD Systems With Boundary Queues Using Colored Petri Nets," *IEEE Transactions on Intelligent Transportation Systems,* vol. 23, no. 8, pp. 12977-12999, 2022, doi: 10.1109/TITS.2021.3119017.

[33] C. Chen, N. Geroliminis, and R. Zhong, "An Iterative Adaptive Dynamic Programming Approach for Macroscopic Fundamental Diagram-Based Perimeter Control and Route Guidance," *Transportation Science,* 2024.

[34] Q. Guo and X. Ban, "Macroscopic fundamental diagram based perimeter control considering dynamic user equilibrium," *Transportation Research Part B: Methodological,* vol. 136, pp. 87-109, 2020, doi: 10.1016/j.trb.2020.03.004.

[35] T. Lei, Z. Hou, and Y. Ren, "Data-Driven Model Free Adaptive Perimeter Control for Multi-Region Urban Traffic Networks With Route Choice," *Ieee Transactions on Intelligent Transportation Systems,* vol. 21, no. 7, pp. 2894-2905, Jul 2020, doi: 10.1109/tits.2019.2921381.

[36] P. Varaiya, "Max pressure control of a network of signalized intersections," *Transportation Research Part C: Emerging Technologies,* vol. 36, pp. 177-195, 2013.

[37] N. Geroliminis and J. Sun, "Properties of a well-defined macroscopic fundamental diagram for urban traffic," *Transportation Research Part B: Methodological,* vol. 45, no. 3, pp. 605-617, 2011.

[38] X. Peng and H. Wang, "Capturing Spatial-Temporal Traffic Patterns: A Dynamic Partitioning Strategy for Heterogeneous Traffic Networks," *IEEE Access,* 2024.